\documentstyle[12pt]{article}
\title{A Model Existence Theorem for Infinitary Formulas in Metric Spaces}
\author{Carlos E. Ortiz \\ Beaver College}
\begin{document}
\maketitle

\newtheorem{lemma}{LEMMA}[section]
\newtheorem{theo}[lemma]{THEOREM}
\newtheorem{df}[lemma]{DEFINITION}
\newtheorem{example}[lemma]{EXAMPLE}
\newtheorem{cor}[lemma]{COROLLARY}
\newtheorem{remark}[lemma]{REMARK}
\newtheorem{prop}[lemma]{PROPOSITION}

\newcommand{\reals}{\mbox{$\Re$\/ }}
\newcommand{\findef}{\mbox{$\Box$} \vspace{.15in}}
\newcommand{\apptruth}{\mbox{$\models_{AP}$ }}
\newcommand{\noapptruth}{\mbox{$\not\models_{AP}$}}
\newcommand{\uequiv}{\mbox{$\sim _{\cal U} $}}
\newcommand{\concati}{\mbox{$
(\vec v_{1},\vec v_{2}..\vec v_{i}..) $}}
\newcommand{\concatf}{\mbox{$ (f_{1},f_{2},..f_{i}..)$}}
\newcommand{\concatg}{\mbox{$ (g_{1},g_{2},..g_{i}..)$}}
\newcommand{\concatgk}{\mbox{$ (g_{1},g_{2},..g_{k}..)$}}
\newcommand{\concatgki}{\mbox{$ (g_{1}(i),g_{2}(i),..g_{k}(i)..)$}}
\newcommand{\concatfki}{\mbox{$ (f_{1}^{k}(i),g_{2}^{k}(i),..g_{i}^{k}(i)
..)$}}
\newcommand{\concatgj}{\mbox{$ (g_{1},g_{2},..g_{j}..)$}}
\newcommand{\concatfj}{\mbox{$ (f_{1},f_{2},..f_{j}..)$}}
\newcommand{\concatgi}{\mbox{$ (g_{1}(i),g_{2}(i),..g_{i}(i)..)$}}
\newcommand{\concatgji}{\mbox{$ (g_{1}(i),g_{2}(i),..g_{j}(i)..)$}}
\newcommand{\concath}{\mbox{$ (h_{1},h_{2},..h_{i}..)$}}
\newcommand{\concathk}{\mbox{$ (h_{1}(k),h_{2}(k),..h_{i}(k)..)$}}
\newcommand{\concathi}{\mbox{$ (h_{1}(i),h_{2}(i),..h_{j}(i)..)$}}
\newcommand{\concatheq}{\mbox{$ (h_{1}(e_{q}),h_{2}(e_{q}),..h_{i}
(e_{q})..)$}}
\newcommand{\formu}{\mbox{$\Psi$}}
\newcommand{\fee}{\mbox{$\phi$}}
\newcommand{\bbar}{\mbox{$\vec b$}}
\newcommand{\structu}{\mbox{\boldmath $U$\/}}
\newcommand{\xbar}{\mbox{$\vec x$}}
\newcommand{\structa}{\mbox{\boldmath $E $ \/}}
\newcommand{\structone}{\mbox{\boldmath $E_{1}$ \/}}
\newcommand{\structtwo }{\mbox{\boldmath $E_{2}\/$}}
\newcommand{\structb}{\mbox{\boldmath $B$\/}}
\newcommand{\structc}{\mbox{\boldmath $C$\/}}
\newcommand{\structz}{\mbox{\boldmath $Z$\/}}
\newcommand{\structi}{\mbox{\boldmath $E_{i}$\/}}
\newcommand{\structj}{\mbox{\boldmath $E_{j}$\/}}
\newcommand{\structf}{\mbox{\boldmath $F$\/}}

\newcommand{\structbeta}{\mbox{\boldmath $E_{\beta}$\/}}
\newcommand{\structgamma}{\mbox{\boldmath $E_{\gamma}$\/}}
\newcommand{\structn}{\mbox{\boldmath $E_{n}$\/}}
\newcommand{\structk}{\mbox{\boldmath $E_{k}$\/}}
\newcommand{\structp}{\mbox{\boldmath $E_{p}$\/}}
\newcommand{\structsucc}{\mbox{\boldmath $E_{k+1}$\/}}
\newcommand{\structni}{\mbox{\boldmath $E_{n_{i}}\/$}}
\newcommand{\structnk}{\mbox{\boldmath $E_{n_{k}}$\/}}
\newcommand{\structm}{\mbox{\boldmath $E_{m}$\/}}
\newcommand{\structsu}{\mbox{\boldmath $E[S,\cal U ]$\/}}
\newcommand{\structspu}{\mbox{\boldmath $E[S',\cal U ]$\/}}
\newcommand{\structsju}{\mbox{\boldmath $E[S_{j},\cal U ]$\/}}
\newcommand{\ultrastructa}{\mbox{$\prod_{\cal U}$\boldmath $E_{i}$\/ }}
\newcommand{\ultrastructni}{\mbox{$\prod_{\cal U}$\boldmath $E_{ni}$\/}}
\newcommand{\fullstructa}{\mbox{\boldmath $E$ $=
(X,d,F,R,C,F_{\reals})$}}
\newcommand{\fullstructb}
{\mbox{\boldmath $B$ $=(Y,\mbox{$d,F_{2}, R_{2},C_{2}, (F_{\reals})_{2}$)}}
\newcommand{\fullstructsu}
{\mbox{\boldmath $E$ $[S,\cal U ]$}=(S/\equivu,$$d,
F, R, C, F_{\reals} $)}}
\newcommand{\fullstructi}
{\mbox{\boldmath $E_{i}$ $=(X^{i},d^{i},F^{i} 
,R^{i},C^{i}, F_{\reals}^{i}$)}}
\newcommand{\fullPhi}{\mbox{$L=(\cal {F}$$, 
\cal {R}$,$\cal C$,$\cal {P}$$, \cal F$$_{\reals})$ }}
\newcommand{\fullultrastructa}
{\mbox{\boldmath $E_{i}$$=((\prod_{i=1}^{\omega} X_{i})/ \uequiv ,
rho,F ,R , C, F_{\reals}$)}}
\newcommand{\APhi}{\mbox{$ \Phi^{app}$}}
\newcommand{\verify}{\mbox{$\models$ }}
\newcommand{\almostsub}{\mbox{$\subseteq^{*}$}}
\newcommand{\strictalmostsub}{\mbox{$\subset^{*}$}}
\newcommand{\orderka}{\mbox{$\unlhd_{K}$}}
\newcommand{\qed}{\mbox{\rule{2mm}{2mm} \vspace{.10in}}}
\newcommand{\concathj}{\mbox{$ (h_{1},h_{2},..h_{j}..)$}}
\newcommand{\nconcathji}{\mbox{$ (h_{1}(i),h_{2}(i),..h_{j}(i)..)$}}
\newcommand{\concatb}{\mbox{$ (\vec b_{1},\vec b_{2},..\vec b_{i}..)$}}
\newcommand{\concatj}{\mbox{$ (\vec v_{1},\vec v_{2},..\vec v_{j}..)$}}
\newcommand{\concatn}{\mbox{$ (\vec v_{1},..\vec v_{n})$}}

\begin{abstract}
We prove a Model Existence Theorem for a fully  infinitary logic $L_{A}$
for metric structures. This result is based on a generalization of
 the notions of approximate formulas and approximate truth in 
 normed structures introduced by Henson (\cite{Henson1}) and studied 
  in different forms by Anderson 
(\cite{Anderson}) and Fajardo and Keisler (\cite{Fajardo&Keisler1}).
This theorem extends Henson's Compactness Theorem 
for approximate truth in normed space structures to infinitary formulas.

\end{abstract}

\section{Introduction}

In 1976 Henson in \cite{Henson} introduced a logic
of positive bounded formulas in Banach spaces to study the relationship
between a Banach space $E$ and its nonstandard hull $H(E)$. This logic
$L_{PB}$ is based on a first order language $L$ containing a binary
function symbol +,  unary predicate symbols $P$ and $Q$ to be interpreted
as the closed unit ball and the closure of its complement and, for every
rational number $r$, a unary function symbol $f_{r}$ to be interpreted
as the operation of scalar multiplication by $r$.
$L_{PB}$ is closed under finite conjunction, disjunction and
bounded quantification of the form ($\exists x)(P(x)\wedge$\ldots)
or ($\forall x)(P(x) \Rightarrow$\ldots).

  For any formula $\phi$ in $L_{PB}$ and for every natural number
 $n$ it is possible to define in a purely syntactical way a formula
  $\phi_{n}$ in $L_{PB}$, called the {\em n-approximation \/} of $\phi$.
  Intuitively, $\phi_{n}$ is the formula that results from metrically
  weakening the predicates that appear in $\phi$ in such a way that as
  $n$ tends to
  $\infty$, $\phi_{n}$ approaches $\phi$.  From this notion of
  approximations of formulas  follows the definition of approximate
truth: a formula $\phi$ $\in L_{PB}$ is approximately true in a Banach
space $E$ (denoted by $E \apptruth \phi$) iff for every integer
$n$,  $E \models \phi_{n}$. It can be shown that $\apptruth$ is a weaker 
notion than $\verify$.

   The above definition of $\apptruth$ is the starting point of
   the model theory of Banach spaces. 
The cornerstone of the development of this model theory  
 is Henson's Compactness Theorem (\cite{Henson&Iovino}).
This theorem states that for particular collections of normed spaces models 
(we will call these collections ``classes of models''), the following holds:
 \begin{quote}
\em  for any $\phi \in
L_{PB}$, if every approximation of $\phi$ holds in some model
of the class, then there is a model in the class where $\phi$ 
 is approximately true and is also true. 
\end{quote}

Henson, Heinrich and Iovino developed the model theory of $\apptruth$ in normed spaces 
(\cite{ Henson1,Henson2,Henson3,Heinrich&Henson,Henson&Iovino}). It
 follows from their work 
that $\apptruth$ is the natural semantic notion 
to do model theory in 
 normed structures since many of the usual theorems of Model Theory 
 (Compactness, Lowenheim-Skolem, etc.) hold for normed structures 
 (using the semantic notion of \apptruth). This is not the case 
 for the usual semantic notion of $\verify$ in normed structures 
 (see \cite{Iovino}). Furthermore, $\apptruth$ is a concept that
appears naturally in analysis. For example, ``The continuous map
 $T:E\longmapsto E$ has a fixed point''
   is approximately true in $E$ iff $T$ has an almost fixed point
 (\cite{Aksoy&Khamsi}). 

 However, there is a fundamental limitation to the above approach.
 Many of the relevant properties in functional analysis require 
a full infinitary logic to be expressed 
(e.g. reflexivity, superreflexivity).  In order to extend 
the applicability of  
the tools and ideas of model theory to problems in functional analysis,
it is desirable then to obtain model theoretic results for
 infinitary formulas 
in normed structures. A natural starting point is to obtain
 a Model Existence Theorem that generalizes Henson's Compactness
 Theorem to an infinitary logic and to general metric structures.

Furthermore, in areas like the geometry of Banach spaces, it is useful to 
have tools to construct spaces that satisfy complex infinitary properties. 
Famous problems solved recently, like the distortion problem
 (\cite{Odell&Schlumprecht1}) or
the unconditional basis problem
 (\cite{Gowers&Maurey}), require the construction of pathological Banach
 spaces.  Problems still open, like the Fixed Point Property
 for reflexive spaces, would require (for a negative solution) the
 construction of pathological Banach spaces 
 with infinitary properties. Any model-theoretic tool used to construct
 metric spaces  satisfying complex infinitary properties should be based
 on a Model Existence Theorem for infinitary logic in metric spaces.

In this paper we extend Henson's Compactness Theorem  to a full infinitary
 logic ($L_{A}$)
 for metric structures.  $L_{A}$ is closed under countable conjunctions, 
negation and bounded existential quantification over countable many 
variables. This logic was introduced in \cite{Ortiz} 
and studied extensively in \cite{Ortiz1}. The main tools  used to 
obtain the Model Existence Theorem are  the notions of
approximate formulas and approximate truth for  the formulas in $L_{A}$
 and the concept of a uniform collection of models.

 As mentioned above, the logic $L_{A}$ as well as the notion of
 approximate truth for $L_{A}$ were introduced in \cite{Ortiz}.
 For the sake of completeness, we define them here. The logic $L_{A}$
 is defined in Section~\ref{signature Phi}. In this section we also
 define the analog of a collection of models in the classical first
 order logic: the uniform collection of models. In 
 Section~\ref{approximate formulas} we introduce a notion of  
 approximate truth for the formulas in $L_{A}$ that extends Henson's. 
The key idea here is to replace Henson's notion  
 of a sequence of approximate formulas $\{\phi_{n}:n \in \omega\}$ 
 (for formulas $\phi \in L_{PB}$) by a tree of approximate formulas 
 $\{\phi_{h,n}: h \in I(\phi) \mbox{ } n \in \omega\}$ where each 
 sequence $\phi_{h,1},\phi_{h,2},\ldots\phi_{h,n},\ldots$ is a branch of  the
 tree. In this way a sentence $\phi$ is 
 approximately true in a metric structure $E$ iff there exists a branch 
 of the tree of approximations of $\phi$ such that all the finitary formulas 
 in this branch hold in $E$. For formulas in $L_{PB}$ it can be shown that the branches of the tree of approximations are all logically equivalent. In this case, our notion of approximate truth coincides with Henson's original one.

 In Section~\ref{comply} we 
 prove a Model Existence Theorem for $L_{A}$ along the following 
 lines:
 \begin{quote}
 If there exists a branch $\phi_{h,1},\ldots\phi_{h,n},\ldots$ of the tree
of approximate formulas of $\phi$ such that for every $n$ there 
is a model $\structn$  in the uniform collection of models $\cal M$ with
$\structn \verify \phi_{h,n}$ then there exists a model \structa 
in $\cal M$ such that : $\structa \verify \phi$ and $\structa \apptruth 
\phi$.
\end{quote}
 
\section{The Logic $L_{A}$ \label{signature Phi} }

We will use the symbol $\reals$ to refer to the metric space of the real 
numbers with the usual metric. Likewise, the symbol $Q$ denotes 
the rationals, and  $Q^{+}$ the rationals strictly greater than zero.

We will use $\structa, \structb$ to denote metric models, letters 
$E,F, \ldots$ to denote first order models and 
letters $A,X, Y,\ldots$ to refer to sets. We reserve 
the letters $f,g$ for function (or function symbols), 
and $\epsilon, \delta$ for real numbers. Likewise, we 
use $m, n, r,\ldots$ for integers.

We can 
visualize a metric language as being made of two pieces:
 a classical first order model and  a metric part.

\begin{df} Definition of a Metric Language
\end{df}

  Fix a first order language $L$ made of function, relation
 and constant symbols (with finite arity). Fix also $M_{\Re}$
 a collection of function symbols disjoint from the function
 symbols in $L$. Let $D_{\Re}$ be the collection 
of all the compact sets in $\Re$ with the usual metric.

    $L(M_{\Re}, D_{\Re})$, the metric language induced by $L$,
 $M_{\Re}$ and $D_{\Re}$, is the collection of all the relation,
 function and constants symbols appearing in $L \cup M_{\Re}$, as
 well as the collection of all the compact sets in $\Re$.

The function symbols in $M_{\Re}$ are called the real valued function 
  symbols of the language. These function symbols are going to be 
 interpreted  as maps taking values in \reals. There is a fixed 
 real valued function symbol $\rho \in M_{\Re}$ 
with arity two. This symbol will be interpreted in every metric model as
 a metric function.\findef

\begin{example} $\empty$
\label{example1}\end{example}

A natural metric language to deal with normed structures is 
$$L(\{\rho(x,y), ||x||\}, D_{\Re})$$ 

Here $L$ is the first order language 
$L=(\{(x+y)\} \cup \{ r(x):
r \in Q^{+}\}, \{ B_{q}: q\in Q \},\{0\})$
with $x+y$ a binary model valued function symbol, 
 $r(x)$ (for 
$r \in Q$) a unary model valued function symbol,  $0$ a constant 
and 
$B_{q}$ (for every positive rational $q$) a unary predicate (to be 
interpreted in the normed structures as the ball of radius $q$ 
centered at the 
origin).

 The metric part of the logic is made of all the compact 
subsets of \reals  (the set $D_{\Re}$), of the binary real valued
 function symbol $\rho(x,y)$  and of 
the unary real valued function symbol $||.||$.

We will call this metric language the language for normed structures. 
Note that every normed space  can be seen in a 
natural way as a structure in this language.\findef

We remark that this definition of metric language can be 
extended naturally to include, for any  complete metric 
space $(X,\rho)$, $(X,\rho)$-valued function symbols, and 
the collection of the compact metric subsets of $(X,\rho)$. 
All the results proved in this paper extend naturally to this 
setting (see \cite{Ortiz1}). For simplicity's sake, we refrain 
of doing so here.

Note that the collection of the compact subsets of the reals ($D_{\Re}$) 
appears in every metric language. In order to reduce 
notation clutter we  omit to mention 
this set   
in the descriptions of the metric languages. It is understood 
that every metric language 
contains this collection.

For the rest of this section we fix a metric language $L(M_{\Re})$.

We define now a logic $L_{A}$ over the metric language $L(M_{\Re})$.
 The terms and the values of the terms are defined by induction.

\begin{df} Definition of terms in $L_{A}$
\end{df}

\begin{enumerate}
\item Let $\{x_{i}| i \in \omega_{1}\}$ be a fixed collection
of variables. The variables $x_i$ are terms. Their arity is 1.
They are considered model valued terms.
 \item Given a collection of model valued terms \{$t_{i}$: $i\leq n$\}
  and a real valued function symbol
  (or a model valued function symbol)
 $f$ with arity $n$, $f(\vec t)=
 f(t_{1},\ldots t_{n})$
 is a real valued term (or model valued term) with arity 
 given by the cardinality of the collection
 of variables in $\{t_{i}:i \leq n\}$. \findef
\end{enumerate}

  In summary, terms have finite arity and they can be model
 valued or real valued. 

We are ready to define the formulas of the logic $L_{A}$
by induction.  As usual, for every formula $\phi$, 
we will use the notation $\phi(\vec x)$ to 
mean that the free variables of $\phi$ are among the components
of the vector $\vec x$. Likewise, $\phi(\vec x_{1},\vec
x_{2},\ldots)$ means that the free variables of $\phi$ are among the
components of the vectors $\vec x_{1},\vec x_{2},\ldots$.

 \begin{df} Atomic formulas in $L_{A}$
 \end{df}

 An atomic formula in $L_{A}$ would be any expression of one of the 
following forms:
 \begin{enumerate}
 \item  $C(\vec t$) ,
where $C$ is a relation symbol in $L$ and $\vec t$ is a vector of model 
valued terms such that
its arity  agrees with the arity of $C$.
\item $B( t)$, where $B \in D_{\Re}$ is a 
compact subset of the reals and $t$ is 
a real valued term. \findef
\end{enumerate}

 \begin{df} Definition of formulas in $L_{A}$
 \end{df}
  \begin{enumerate}
   \item An atomic formula is a formula in $L_{A}$.

    \item If $\phi_1,\phi_2,\ldots\phi_i$\ldots ($i{<}\omega$)
     is a collection of formulas in $L_{A}$, then for every
     $n$,
      \begin{displaymath}
      \bigwedge_{i=1}^{n} \phi_{i} \mbox{ and } \bigwedge_{i=1}^{\infty}
      \phi_{i}  
      \end{displaymath}
      are in $L_{A}$.
      \item If $\phi$ is a formula in $L_{A}$ then $\neg\phi$ is
      also a formula in $L_{A}$.

\item Consider a  formula $\phi( \vec y_{1},\ldots \vec y_{n},\ldots,
x)$ in $L_{A}$. Let $\vec
K=(K_{1},K_{2},\ldots)$ be a corresponding vector of relation symbols.
 The following formula is in $L_{A}$:
  \begin{displaymath}
    \exists (\vec y_{1},\ldots \vec y_{n},\ldots )\in (K_{1},\ldots 
    K_{n},\ldots) (
      \phi( \vec y_{1},\ldots  \vec y_{n},\ldots ,\vec x )). \mbox{ \findef}
        \end{displaymath}
\end{enumerate}

We will abbreviate $\neg \bigwedge \neg$ by $\bigvee$ and $\neg
(\phi \wedge \neg \psi)$ by $ \phi \Rightarrow \psi$. To avoid very long formulas we will also abbreviate
\begin{displaymath}
\exists ( y_{1}\ldots  y_{n}\ldots ) \in
(K_{1},\ldots K_{n},\ldots )(\phi( x_{1},\ldots  x_{n},\ldots ,\vec y))
\end{displaymath}
by $\exists \vec y\in \vec K \phi(\vec x,\vec y)$. Likewise $\neg
\exists \vec y \in \vec K \neg \phi(\vec y,\vec x)$ will be
abbreviated by $\forall \vec y\in \vec K \phi(\vec
y,\vec x)$.

 Finally, given a countable set $A$ 
 with a fixed enumeration $A=\{a_{1},\ldots,a_{n},\ldots\}$ 
 and  countable formulas $\{\phi_{a}\}_{a\in A}$ we understand 
 by $\bigwedge_{a\in A} \phi_{a}$ the formula 
 $\bigwedge_{n=1}^{\infty} \phi_{a_{n}}$. 
 Likewise, for any  arbitrary integer $r$, 
 we understand by $\bigwedge_{a \in A\uparrow r} 
 \phi_{a}$ the formula $\bigwedge_{n=1}^{r} \phi_{a_{n}}$.

Notice that the formulas in $L_{A}$ admit bounded quantification over
infinitely many variables. This enhances the expressive power of
$L_{A}$, as the next example shows.

\begin{example} Expressive power of $L_{A}$
\end{example}
 {\bf The normed structure \structa  is reflexive. \/}

For any Banach space $(X,||.||)$, let $B_{1}$ denote the unitary ball.
 A characterization of reflexivity due to James (\cite{James})
(see also \cite{Van Dulst}) 
 that does not require any 
mention of the dual is the following:

{\em 
A Banach space $(X,||.||)$  is reflexive 
iff\/} $$ \forall \epsilon > 0 \mbox{ }
\forall \{x_{i}\}_{i=1}^{\infty} \subseteq  B_{1}\mbox{ } \exists k
 \in \omega \mbox{ }
dist[conv(\{x_{1},\ldots,x_{k}\}),conv(\{x_{k+1},\ldots\})] 
\leq \epsilon$$ 

Here, for any set  $A \subset X$,  conv(A) is the convex hull spawned 
by $A$. Similarly,  given two sets $A,B \in X$, $dist[A,B]=inf\{||x-y||
 : x \in A, y \in B\}$.

Recall the language for normed structures discussed in
 Example~\ref{example1}, 
$L=(\{\rho(x,y), ||.||\})$, with $L$ the first order language 
$$L=(\{(x+y)\} \cup \{ r(x):
r \in Q^{+}\}, \{ B_{q}: q\in Q^{+} \},\{0\})$$

 The following sentence of $L_{A}$ expresses reflexivity in this language:
\begin{displaymath}
\bigwedge_{n=1}^{\infty} \forall (x_{1},\ldots x_{m},\ldots )\in (B_{1},
B_{1},\ldots B_{1},\ldots) 
\end{displaymath}
\begin{displaymath}
\bigvee_{k=1}^{\infty} \bigvee_{r=1}^{\infty} 
\bigvee_{\vec a \in CO(k)} 
\bigvee_{\vec b \in CO(r)} 
||\sum_{i=1}^{k} a_{i}x_{i} -\sum_{j=1}^{r} b_{j}x_{k+j}|| 
\leq (1/n)
\end{displaymath}
Here, $\forall s \in \omega$, $CO(s)$ 
is the subset of $Q^{s}$ made of all the $s$-tuples 
$(a_{1},\ldots a_{s})$ such that $\sum_{i=1}^{s} a_{i}=1$ and 
$a_{1},\ldots a_{s} \geq 0$. 

 {\bf The normed structure \structa  is uniformly convex. \/}
A normed space $(X,||.||)$ is uniformly convex iff $\forall \epsilon > 0$ 
 $\exists \delta > 0$ such that for every $x,y $ in the unitary sphere,
$||x-y|| \geq \epsilon \Rightarrow \frac{||x+y||}{2} \leq 1-\delta$.

We use again the language for normed structures $L(\{\rho,||.||\})$. 
Here is how we express uniform convexity by a sentence in $L_{A}$:
\begin{displaymath}
\bigwedge_{n=1}^{\infty} \bigvee_{m=1}^{\infty} 
\forall (x,y) \in (B_{1},B_{1}) [ \mbox{ }||x-y||\geq (1/n) 
\Rightarrow (1/2)||x+y|| \leq (1-1/m) \mbox{ }] \mbox{ }\findef
\end{displaymath}

\subsection{Semantics for $L_{A}$}
\label{semantics of LA}

The definition of the structures of $L_{A}$ is a natural
generalization of Henson's notion of a Normed Space Structure
(see\cite{Henson&Iovino}).

 Any model of a metric language $L(M_{\Re})$ can be 
 imagined as a classical 
 first order model plus interpretations for the metric part of $M_{\Re}$.

\begin{df} Definition of Model
\end{df}
 A model for $L(M_{\Re})$
  is a pair $ \structa=(E,F_{\Re})$
  where:
\begin{enumerate}
\item $E$ is a model of the first order language 
$L$ with interpretations for the function, relation and constant symbols 
in $L$. Let us call $X$ the universe of $E$. 

\item  the real valued function symbol $\rho$ is interpreted as a 
metric in $X$. Furthermore, every interpretation of a model valued 
function symbol $f \in L$ is a continuous function $f^{\structa}$ 
 (with respect to the product topology
induced by $\rho$ in $X^{a_{f}}$) from $X^{a_{f}}$ to $(X, \rho)$.

\item $F_{\Re}=\{ f^{\structa}| f\in M_{\Re}\}$ with 
the property that 
for every  real  valued $ f\in M_{\Re}$ with arity $a_{f}$,
 the interpretation $f^{\structa}$
is a continuous function, with respect to the product topology
induced by $\rho$ in $X^{a_{f}}$, from $X^{a_{f}}$ to $\Re$. \findef
\end{enumerate}

  Once the interpretations for the symbols of $L(M_{\Re})$ 
  are defined
  one can define the interpretations of the terms in the model,
  $t^{\structa}$, in the natural way. When the model and the 
  interpretations of the elements of the language are clear 
  from the context we will drop the symbol
   $t^{\structa}$ and use $t$.

The definition of satisfaction follows in the standard way.
\begin{df} $\empty$
\end{df}
 The truth ($\models$) relation in the model $\structa=(E, F_{\Re})$
  is constructed by induction in formulas in the usual way from
  the truth definition for the atomic case:
\begin{itemize}

\item       Let $K(\vec t(\vec x))$ be an atomic formula in $L_{A}$ with 
$K \in L$. Let $\vec
       a\in X^{|\vec x|}$. $\structa 
       \verify K(t(\vec a))$ 
       iff $E$ $\models$ $K ( \vec t( \vec a ))$ in the usual first 
       order sense.
\item Let $B(t(\vec x))$ be an atomic formula with $B$ a compact 
subset of the reals ($B \in D_{\Re}$).
Let $\vec a \in X^{|\vec x|}$. $\structa 
\verify B(t(\vec a))$ 
iff 
it is true that $t^{\structa}(\vec a) \in B$.\findef
\end{itemize}

\subsection{Uniform Collection of Metric Models}

The analog to the notion of collection of models for a first 
order language is given by the notion of uniform collection of 
models in metric languages.

Intuitively, a uniform collection  of models is a collection of
 all models for $L(M_{\Re}) $ satisfying the same uniform 
 continuity requirements for
 the interpretations of the function symbols and  satisfying the
 same uniform ``bound'' for the interpretations of the predicate
 symbols. Recall that in this section we are dealing with a fix
  metric language $L(M_{\Re})$.

\begin{df}
Uniform Collection  of Models
\label{uniformcollection}
\end{df}
 Fix  the following assignments:
\begin{itemize}
\item For every model valued (real valued) function symbol 
$f$ with arity $r$, for every vector
$\vec K=(K_{1},\ldots K_{r})$ of  predicate symbols in $L$, 
fix a  predicate symbol $K^{(f,\vec K)}$
in $L$ (or a compact predicate $D^{(f,\vec K)}$ in
$D_{\Re}$). It is intended for $K^{(f,\vec K)}$ ($D^{(f,\vec K)}$)  to
contain the set  $f(\vec K)$ in every model of 
the class.
\item For every two relation symbols $Q, A$ with arities $a,b$ in $L$, fix a  
predicate symbol $K^{(Q,A)}$ in $L$. It is 
intended for  $K^{(Q,A)}$  to be a uniform cover of the union of 
the projections of the relations $Q,A$ into
the model.
\item
For every model valued (or real valued) function symbol 
$f$ with arity $r$, for every vector 
$\vec K=(K_{1},\ldots K_{a})$
of predicate symbols in $L$, for every
rational $\epsilon >0$ fix a rational $\delta^{(f,\vec K,\epsilon)}
> 0$.
\end{itemize}

A uniform collection  of models for the above assignments 
 is the collection of all the models $\structa=(E,F_{\Re})$ 
 in $L(M_{\Re})$ that satisfies:
 \begin{enumerate}
  \item {\bf Uniform bound for function symbols\/}.
   For every model valued (real valued) 
function symbol $f$ in $L$, for every vector
     $\vec K=(K_{1},\ldots K_{r})$ of predicate
      symbols in $L$,
       \begin{displaymath}
        \structa \verify \forall \vec x \in 
        \vec K\mbox{  $\empty$}
        K^{(f,\vec K)}(f(\vec x))
         \end{displaymath}
          (or 
\begin{displaymath}
 \structa \verify \forall \vec x \in \vec K \mbox{ $\empty$  }
 D^{(f,\vec K)}(f(\vec x))\mbox{  )}
  \end{displaymath}

  \item {\bf The relation symbols of $L$ form a directed set.\/}
  For every relation symbols $Q,A$ in $L$ with arities $a$, $b$,
  \begin{displaymath}
  \structa \verify \forall \vec x,\vec y \in (Q, A)(
  \bigwedge_{i=1}^{a} K^{(Q,A)}(x_{i}) \wedge
  \bigwedge_{j=1}^{b} K^{(Q,A)}(y_{j}))
  \end{displaymath}
   Furthermore $\structa$ verifies that for
    every $a$ in \structa there exists a  predicate symbol 
    $K$ in $L$ such that $(\structa, F_{\Re}) \verify
    K(a)$.

    \item {\bf Uniform continuity for function symbols\/}.
    For every model valued (or real valued) function symbol 
    $f$ with arity $a$, for every vector  
    $\vec K=(K_{1},\ldots K_{a})$
of predicate symbols in $L$, for every
rational $\epsilon >0$,
\begin{displaymath}
\structa \verify \forall \vec x,\vec y \in (\vec K,\vec K) 
(\bigwedge_{i\leq a} \rho(
x_{i}, y_{i}) < \delta^{(f,\vec K,\epsilon)} \Rightarrow
\rho(f(\vec x),f(\vec y)) \leq \epsilon )
\end{displaymath}
(or
\begin{displaymath}
\structa \verify \forall \vec x,\vec y \in (\vec K,\vec K) 
(\bigwedge_{i\leq a}
|x_{i} - y_{i} |
< \delta^{(f,\vec K,\epsilon)} \Rightarrow
\rho_{Y}(f(\vec x),f(\vec y)) \leq \epsilon )\mbox{   )} \findef
\end{displaymath}
\end{enumerate}

\begin{example} Uniform Class of Normed Space Structures
\label{complete}\end{example}

Consider the metric language 
$L(\{\rho(x,y),||x||\})$ with L the first order language
$L=(\{(x+y) \}\cup \{ r(x):
r \in Q\}, \{ B_{q}: q\in Q^{+}\},\{0\},\}$ and 
with $\{\rho(x,y),||x||\}$ 
as discussed 
in Example~\ref{example1}.

The natural uniform collection of models associated with this 
language and with the normed spaces is spawned by the following 
natural assignments:
\begin{itemize}
\item For $\rho (x,y)$ and for every $B_{q_{1}}$, $B_{q_{2}}$, 
let $$D^{(\rho,(B_{q_{1}},B_{q_{2}}))}=\{ z \in \Re : 
0 \leq z \leq q_{1}+q_{2}\}. $$ 
 For $(x+y)$ and for every $B_{q_{1}}$, $B_{q_{2}}$,
let $$K^{(+,(B_{q_{1}},B_{q_{2}}))}=B_{q_{1}+q_{2}}.$$ 
For $||x||$ and $B_{q}$ let $D^{(||.||,B_{q})}=\{ z \in \Re :
0 \leq z \leq q\}$.

\item 
 For every model valued function $r(x)$, with $r \in Q$, and for every 
$B_{q}$ we assign $K^{(r(.),B_{q})}=B_{|r|q}$. 
For every $B_{q_{1}}$, $B_{q_{2}}$, we assign 
$K^{(B_{q_{1}},B_{q_{2}})}
=B_{q_{1}+q_{2}}$.
\item For $\rho (x,y)$, for every $B_{q_{1}}$, $B_{q_{2}}$ and 
for every rational $\epsilon >0$, let $\delta^{(\rho,(B_{q_{1}},B_{q_{2}})
,\epsilon)}=\epsilon/2$. For $(x+y)$, for every 
$B_{q_{1}}$, $B_{q_{2}}$ and for every rational 
 $\epsilon >0$, 
let $\delta^{(+,(B_{q_{1}},B_{q_{2}}),\epsilon)}=\epsilon/2$. 
 For $||x||$, $B_{q}$ and a rational $\epsilon > 0$, 
let $\delta^{(||.||,B_{q},\epsilon)}=\epsilon$.
  For every function $r(x)$, with $r \in Q$, for every 
$B_{q}$ and every rational $\epsilon > 0$, 
we assign $\delta^{(r(.),B_{q})}=\epsilon/r$.

\end{itemize}

It is clear that every model $\structa=(E, \{\rho, ||.||\})$ of 
$ L(M_{\Re})$
such that  
$E=(X,\{0\},+,\{r(x): r \in Q\})$ is a vector space (over $Q$) and
$||.||$ is a norm in $X$, belongs to the above uniform collection of 
models.\findef

\section{Approximate Formulas in $L_{A}$}
\label{approximate formulas}

We define a notion of approximate formula (for the formulas in
$L_{A}$) that extends Henson's definition
(see~\cite{Henson&Iovino}). 

Our intention is to generate approximations of all the formulas
in $L_{A}$ by using the finitary formulas in $L_{A}$ as building blocks.
Inspired by the Souslin operation used in classical descriptive
set theory (see for example \cite{Kechris}) to generate the
analytic sets from the closed sets in a Polish space, we associate to 
every formula $\phi$ in $L_{A}$ a tree of finitary formulas in 
$L_{A}$ .

Formally, we associate to every formula $\phi$ in $L_{A}$ a set of
indices $I(\phi)$ (the branches of the tree of approximate
formulas) and a set $\{ \phi_{h,n}|
h \in I(\phi),n \in \omega\}$ of finitary formulas in
$L_{A}$. Intuitively, for  every branch $ h\in I(\phi)$,
the approximate formulas in the branch (the collection 
$\{\phi_{h,n}| n \in \omega\}$) are going to
``approach'' $\phi$ as $n$
tends to $\infty$. 

The notions of $I(\phi)$ and $\phi_{h,n}$ were introduced in \cite{Ortiz}
 and studied in detail in \cite{Ortiz1}.

Notation: sometimes, to emphasize that the pair $h,n$ affects
the formula $\phi$, we will write $(\phi)_{h,n}$ for the approximation
$\phi_{h,n}$.

Recall that we are working with a fixed 
metric language $L(M_{\Re})$. Given two formulas $\phi, \sigma$, we 
will write $\phi \equiv \sigma$ if $\phi$ and $\sigma$ are identical formulas.

\begin{df} Definition of approximate formulas in $L_{A}$
\label{defofapprox}
\end{df}

For any formula $\phi(\vec x)$ in $L_{A}$ we define:
\begin{enumerate}
\item  A set of branches $I(\phi)$,
\item  For any $h\in I(\phi)$ and $n\in\omega$, a finitary formula
$\phi_{h,n}(\vec x)\in L_{A}$.
\end{enumerate}

The definition is by induction in formulas as follows:

 {\bf Atomic.\/} For any atomic formula $C(t_{1},\ldots t_{r})$ or 
$D(t)$:
\begin{itemize}
\item The set $I(C(\vec t))=I(D(\vec t))=\{\emptyset\}$
\item If $C \in L$  , then 
for every $h$ in $I(C(\vec t))$, for every integer $n$,
$(C(\vec t))_{h,n}\equiv \exists \vec y \in C 
\bigwedge_{i=1}^{r} \rho (t_{i},y_{i})\leq (1/n)$ (this is 
the $1/n$ metric deformation of $C$).
\item If $B\in \cal M_{\Re}$ is a compact subset of the reals, then the set 
$B_{n}=\{ y \in \reals | \exists x \in B 
\mbox{ with } | x-y| \leq (1/n)\}$ is 
also a compact subset of $\reals$. 
We define then $(B( t))_{h,n}= B_{n}( t)$ 
(this is the (1/n) deformation of the  compact set $B$).
\end{itemize}

{\bf Conjunction.\/} For any countable collection 
$\{\phi_{i}\}_{i=1}^{\infty}$ of
formulas in $L_{A}$, we define:
\begin{itemize}
\item $I(\bigwedge_{i=1}^{\infty} \phi_{i}(\vec
x))=\prod_{i=1}^{\infty}
I(\phi_{i})$ (the cartesian product of the $I(\phi_{i})$ ).
\item For every $h$ in $I(\bigwedge_{i=1}^{\infty} \phi_{i})$,
for every integer $n$,
\begin{displaymath}
(\bigwedge_{i=1}^{\infty} \phi_{i})_{h,n}\equiv
\bigwedge_{i=1}^{n} (\mbox{ }(\phi_{i})_{h(i),n}\mbox{ })
\end{displaymath}
\end{itemize}
{\bf Negation.\/} For any formula $\phi$ in $L_{A}$, we have:
\begin{itemize}
\item $I(\neg \phi) \subseteq (I(\phi)\times \omega)^{\omega}$
is the collection of all maps $f=(f_{1},f_{2})$ with the
following ``weak'' surjectivity property:
\begin{displaymath}
\forall h \in I(\phi) \mbox{ } \exists s \in \omega \mbox{, }
(\phi)_{h,f_{2}(s)} \equiv (\phi)_{f_{1}(s),f_{2}(s)}
\end{displaymath}
\item For every $h \in I(\neg \phi(\vec x))$, for every integer
$n$,$ (\neg \phi)_{h,n} \equiv \bigwedge_{i=1}^{n} \neg
\mbox{ }(\phi_{h_{1}(i),h_{2}(i)})$
\end{itemize}
{\bf Existential.\/} For every formula $\phi(\concati,\vec
x)$, for every vector of bounding sets $\vec
K=(K_{1},K_{2},\ldots K_{i}\ldots )$ of corresponding arities, we have:
\begin{itemize}
\item $I(\exists \vec v \in \vec K \phi(\vec v, \vec x)\mbox{
})=I(\phi(\vec v, \vec x))$
\item For every $h$ in $I(\exists \vec v \in
\vec K \phi(\vec v,\vec x))$, for every
integer $n$, we have that
\begin{displaymath}
(\exists \vec v \in \vec K
\phi(\vec v, \vec x))_{h,n} \equiv
\exists \vec v \in \vec K (\phi(\vec v,\vec x))_{h,n}.
\mbox{ \findef}
\end{displaymath}
\end{itemize}

We  introduce some notation.

We call  $I(\phi)$ {\bf the set of branches\/} of $\phi$. The
formulas $\phi_{h,n}$ are the {\bf approximate formulas\/}  of
$\phi$, and the collection $T(\phi)=\{\phi_{h,n}:h\in I(\phi) 
\mbox{, } n \in \omega\}$
is the {\bf  approximation tree\/} of the formula $\phi$. 

\begin{df} Approximate Truth
\end{df}

Fix a model $\structa$. Let $\phi(\vec x)$ be an arbitrary 
formula in $L_{A}$. Let $\vec b$ be a vector of
elements of $X$.
 We say that
 \begin{displaymath}
 \structa \apptruth \phi ( \bbar ) \mbox{( \structa approximately satisfy
 $\phi$) }
 \end{displaymath} 
iff $\exists h\in I(\phi(\vec x)) \mbox{ } \forall 
n\in\omega \mbox{, }
\structa \verify
\phi_{h ,n } ( \bbar )$. 

Equivalently, $\structa \apptruth \phi (\vec b)$ 
iff there is a branch $h$ of the approximation tree of  $\phi(\vec x)$
 such that for every integer n, 
$\structa \verify \phi_{h,n}$.

We say that $ \structa \apptruth^{h} \phi(\bbar)$ iff 
for every integer $n$, $\structa \apptruth \phi_{h,n}(\bbar)$.\findef

We illustrate these definitions with some examples.

\begin{example} Isometric Inclusion
\label{representability}\end{example}

Let $(Y,|.|)$ be a separable normed space generated by the
countable set of independent vectors
$\{e_{1},e_{2},\ldots e_{n},\ldots \}$ of norm 1. It is easy to see that
$(Y,|.|)$ is isometrically embedded into a Banach space 
$(X,||.||)$ iff there exist unitary vectors 
$v_{1},v_{2}\ldots,v_{n},\ldots$ in $X$ such that:
 \begin{displaymath}
\bigwedge_{m=1}^{\infty} \bigwedge_{\vec a \in Q^{m}}
| \mbox{ }||\sum_{i=1}^{m}
a_{i} v_{i}|| - |\sum_{i=1}^{m} a_{i}e_{i}|\mbox{ }| \leq 0
\end{displaymath}

Consider the language for normed structures $L(\{\rho(x,y),||x||\})$
 discussed in Example~\ref{example1}.
Recall that any normed space can be seen as a structure of 
this language with the natural interpretations.
By $\vec B_{1}$ we denote the vector:
$(B_{1},B_{1},\ldots )$.

We can express that $(Y,|.|)$ is isometrically embedded in  a model
 \structa by the following sentence in $L_{A}$:
\begin{displaymath}
\structa \verify \exists \vec x \in \vec B_{1} 
\bigwedge_{m=1}^{\infty} \bigwedge_{\vec a \in Q^{m}}
| \mbox{ }||\sum_{i=1}^{m}
a_{i} x_{i}|| - |\sum_{i=1}^{m} a_{i}e_{i}|\mbox{ }| \leq 0
\end{displaymath}

On the other hand, the statement
\begin{displaymath}
\structa \apptruth \exists \vec x \in \vec B_{1}
\bigwedge_{m=1}^{\infty} \bigwedge_{\vec a \in Q^{m}}
|\mbox{ }||\sum_{i=1}^{m}
 a_{i} x_{i}|| - |\sum_{i=1}^{m} a_{i}e_{i}|\mbox{ }| \leq 0
  \end{displaymath}
    is equivalent to the statement (since the set of branches
 $I(\phi)=\{\emptyset\}$) that
    for every integer $r$
    \begin{displaymath}
    \structa \verify \exists \vec x \in \vec B_{1}
    \bigwedge_{m=1}^{r} \bigwedge_{\vec a \in
    Q^{m}\uparrow r}
    | \mbox{  }||\sum_{i=1}^{m} a_{i}x_{i}|| -|\sum_{i=1}^{m} a_{i}
    e_{i}| \mbox{  } | \leq 1/r
    \end{displaymath}
It follows easily that this, in turn, is equivalent to the fact 
that for every finite dimensional subspace $G$ of $Y$ and for every
$\epsilon > 0$ there is a finite dimensional subspace $H$ in \structa
that is $\epsilon$-isomorphic to $G$. This is exactly the
definition of the concept of finite representability of a normed
space $Y$ into $X$ (see \cite{Lindenstrauss&Tzafriri}).

In summary: the approximate validity of the statement ``$Y$ is
isometrically embeddable in \structa'' is equivalent to the
statement that $Y$ is finitely representable in \structa. It is
well known that the two notions are not equivalent (see for
example \cite{Lindenstrauss&Tzafriri}). This shows that the
concept of approximate truth is different from validity. \findef

\begin{example} $\empty$
\label{firstexample} \end{example}

Consider a  formula $\phi $ such that 
 $I(\phi)=\{o\}$ has cardinality one. 
 By definition  $I(\neg \phi)=
(\{o\}\times \omega)^{\omega}$. Clearly 
this set can be identified with $(\omega)^{\omega}$. 
Using the definition of approximate truth, we see 
that for every $h \in (\omega)^{\omega}$, for every 
integer $n$, $(\neg \phi)_{h, n}= \bigwedge_{s=1}^{n} 
\neg \mbox{ } (\phi)_{o, h(s)}$.

It follows that for  every structure \structa, $\structa 
\apptruth \neg \phi$ iff $\exists 
h \in (\omega)^{\omega}$ such that for every integer $n$, 
$\structa \verify \bigwedge_{s=1}^{n} \neg 
\mbox{ } (\phi)_{o, h(s)}$. Clearly this is equivalent to 
say that there exists an integer $r$ such that,
$\structa \verify \neg  \mbox{ } (\phi)_{o,r}$. \findef

\begin{example} $\empty$
\end{example}

Consider a countable collection of formulas $\{\phi^{i} : i \in A\}$ 
such that for every $i$, 
 $I(\phi_{i})=\{o_{i}\}$ has cardinality one. By definition, 
 using Example~\ref{firstexample},  $I(\bigwedge_{i \in A} \neg  \phi^{i})$ 
 can be identified with $(\omega)^{A \times \omega}$. Using the
 definition of approximate truth, we see that for every $H \in
 (\omega)^{A \times \omega}$, for every integer $n$,
 $(\bigwedge_{i \in A} \neg 
\phi^{i})_{H,n}= \bigwedge_{i \in A \uparrow  n} \bigwedge_{s=1}^{n} 
\neg \mbox{ } (\phi^{i})_{o_{i},H(i,s)}$.

It follows that for every structure \structa, $\structa \apptruth 
\bigwedge_{i \in A} \neg \phi^{i}$ iff 
there exists $H: A \times \omega \rightarrow \omega$ such that for 
every integer $n$, \structa 
$\verify \bigwedge_{i\in A \uparrow n } \bigwedge_{s=1}^{n} \neg 
\mbox{ } (\phi^{i})_{o_{i},H(i,s)}$. Clearly again, this is equivalent
 to say that for every $i \in A$ there exists an integer $n_{i}$ such that 
$\structa \verify \neg \mbox{ } (\phi^{i})_{o_{i}, n_{i}}$. \findef

\begin{remark} $\empty$
\label{logichenson}
\end{remark}

It is not difficult to see by induction in formulas that if the
 formula $\phi$ in $L_{A}$ is a positive formula
(i.e. a formula without negation) then $I(\phi)=\{ \emptyset \}$
has only one element, i.e. there is only one branch in the
approximation tree of $\phi$.  
For formulas that are not positive, the
set of paths of the approximation tree has cardinality bigger than one
($2^{\omega}$).

Consider the collection $L_{PBA} \subseteq L_{A}$ of all the formulas containing the atomic formulas and closed under countable conjunction, finite disjunction, existential and universal bounded quantification. This collection clearly contains the Henson's logic $L_{PB}$. Furthermore, it can be shown directly (see~\cite{Ortiz1}) that the branches of the tree of approximation  for the formulas $\phi \in L_{PB}$ are all logically equivalent in the following sense: $\forall h,g \in I(\phi)$ $\structa \apptruth ^{h} \phi $ iff $\structa \apptruth ^{g} \phi$. This means that the tree of approximations essentially  contains only one branch $h$.  It is not hard then to verify that for every formula in $L_{PB}$ our  notion of \apptruth  and Henson's notion of approximate truth coincide.  \findef

In the next lemma, we prove two basic properties of the tree of 
approximations of any formula in $L_{A}$. In particular we show that 
the number of approximate formulas $\phi_{h,n}$ for any formula $\phi$ 
 is countable.
We remark that this lemma was already proved in \cite{Ortiz}. We include the proof here for the sake of completeness.

\begin{lemma} $\empty$
\label{setDense(phi)}

 For every formula $\phi(\vec x) \in L_{A}$:
  \begin{enumerate}
 \item
   For any model \structa, for any vector of elements in $\vec b$ in $\structa$, for any integer $n$  and
   for all $h\in I(\phi)$,
\begin{displaymath}
\structa\models \phi_{h,n+1}(\vec b)\Rightarrow
\phi_{h,n}(\vec b)
\end{displaymath}
\item For every formula $\phi$ there exists a set $Dense(\phi)
\subseteq I(\phi)$ at most countable with the following property:

For every branch $h$ in $I(\phi)$, for every
integer $n$, there exists a branch $g$ in $Dense(\phi)$ such that:
$\phi_{h,n}(\vec x) \equiv \phi_{g,n}(\vec x)$
\findef
\end{enumerate}
 \end{lemma}

  PROOF: 

1) By induction on the complexity of the formulas and
   Definition~\ref{defofapprox} of approximate formulas. Left to
   the reader.

   2) By induction on the complexity of the formulas in $L_{A}$.
    
{\bf Atomic Formulas}. Let $C(\vec t)$ (or $B(t)$) be an
atomic formula.  Then we know that
$I(C(\vec t))=\lbrace \emptyset\rbrace$  (or $I(B(t))=\{\emptyset \}$) and
$\forall n\in\omega\forall h\in I(C(\vec t))$ (or in $B(t)$):
\begin{displaymath}
(C(\vec t))_{h,n}=C_{n}(\vec t) \mbox{ (or $(B(t))_{h,n}=B_{n}(t))$ }
\end{displaymath}

Define then $Dense(C(\vec t))=I(C(\vec t)) = \{\emptyset\}$ (or $Dense(B(t))=I(B(t)) =\{\emptyset\}$). It is easy to see that
this set verifies the desired property.

 For the connectives and quantifier steps let us assume as
 induction hypothesis that for formulas $\psi$ of less complexity
 than the formula $\phi$, $I(\psi)$ and $Dense(\psi)$ satisfy the
 desired properties.

 {\bf Conjunction\/}. Consider the formula $ \phi(\vec x)=
 \bigwedge_{i=1}^{\infty} \phi_i(\vec x)$.

Recall that $I(\phi)=\prod_{i=1}^{\infty} I(\phi_i) $
 and that for all integers $n$ and for all ${h}\in I(\phi)$, $ (\phi(\vec x))_{h,n}=\bigwedge_{i=1}^{n}\mbox{ } (\phi_i(\vec
   x))_{h(i),n}$.

    The construction of $Dense(\bigwedge_{i=1}^{\infty} \phi_{i})$ is as
    follows:

    $\forall n \in \omega$, let
    \begin{displaymath}
    p_{n}: \prod_{i=1}^{\infty} Dense(\phi_{i}) \longmapsto
    \prod_{i=1}^{n}
    Dense(\phi_{i})
    \end{displaymath}
    be the usual projection map. By the induction hypothesis, $\forall i
\in \omega$, $Dense(\phi_{i})$ is at most countable. It is possible
then to find, for every integer $n$, countable sets $Q_{n}
\subseteq \prod_{i=1}^{\infty} Dense(\phi_{i}) \subset
I(\bigwedge_{i=1}^{\infty} \phi_{i})$ such that:
 \begin{displaymath}
  p_{n}(Q_{n})= \prod_{i=1}^{n} Dense(\phi_{i})
   \end{displaymath}
   We define $Dense(\bigwedge_{i=1}^{\infty} \phi_{i})=\cup_{n=1}^{\infty}
   Q_{n}$.

   Clearly $Dense(\bigwedge_{i=1}^{\infty} \phi) \subseteq
   I(\bigwedge_{i=1}^{\infty} \phi)$ and is at most countable.
   Furthermore, by the induction hypothesis, for every
   $h=(h_{1},\ldots h_{i},\ldots ) \in I(\bigwedge_{i=1}^{\infty}
   \phi_{i})$, for every integer $n$, there exists $(g_{1},\ldots g_{n})
   \in \prod_{i=1}^{n} Dense(\phi_{i})$ such that:
   \begin{displaymath}
   \bigwedge_{i=1}^{n} (\phi_{i})_{h_{i},n}\equiv
   \bigwedge_{i=1}^{n} (\phi_{i})_{g_{i},n}
   \end{displaymath}
It follows that for every integer $n$ there exists $g \in
Dense(\bigwedge_{i=1}^{\infty}
\phi_{i})$ such that:
\begin{displaymath}
(\bigwedge_{i=1}^{\infty} \phi_{i})_{h,n}\equiv
\bigwedge_{i=1}^{n}
(\phi_{i})_{h_{i},n} \equiv \bigwedge_{i=1}^{n}
(\phi_{i})_{g_{i},n}
\equiv (\bigwedge_{i=1}^{\infty} \phi_{i})_{g,n}
\end{displaymath}
This is the desired result.

{\bf Negation}. Consider the formula
$\phi(\vec x)$= $\neg\psi(\vec x)$ in $L_{A}$.

Recall that $I(\neg \psi) \subseteq (I(\psi)\times
\omega)^{\omega}$ is the collection of all the functions
$f=(f_{1},f_{2})$ with the property that $ \forall h \in I(\psi) \exists s $,   $\psi_{h,f_{2}(s)}
\equiv \psi_{f_{1}(s),f_{2}(s)}$.

Recall also  that  $\forall n\in \omega \mbox{ }\forall
h=(h_1,h_2)\in I(\neg\psi)$ $
(\neg\psi)_{h,n}=\bigwedge_{i=1}^{n}
\neg (\psi_{h_1(i),h_2(i)})$.

The construction of the set $Dense(\neg \psi)$ is as follows:

Let $A \subseteq I(\neg \psi)$ be the collection of all the
functions $f=(f_{1},f_{2}) \in I(\neg \psi)$ such that
$Image(f_{1})\subseteq Dense(\psi)$. From the induction hypothesis it
follows that $\forall n \in \omega$ the set
$(Dense(\psi)\times \omega)^{\{1,\ldots n\}}$ is at most countable. For
any integer let
\begin{displaymath}
proj_{n}:A \longmapsto
(Dense(\psi)\times \omega)^{\{1,\ldots n\}}
\end{displaymath}
be the natural projection map. It is possible to find for
every integer $n$ a countable set $O(n) \subseteq A \subseteq
I(\neg \psi)
\subset (I(\psi)\times \omega)^{\omega}$ such that
\begin{displaymath}
proj_{n}(O(n))=proj_{n}(A)
\end{displaymath}
Let then $Dense(\neg \psi)=\cup_{n=1}^{\infty} O(n)$.

Let us verify the desired property. First, it is clear that
$Dense(\neg \psi)$ is at most countable. Fix now $ h=(h_{1},h_{2})
\in I(\neg \psi)$ and an integer $n$.

By induction hypothesis for every integer $s$ there exists $g_{s}
\in Dense(\psi)$ such that
\begin{displaymath}
\psi_{h_{1}(s),h_{2}(s)}\equiv \psi_{g_{s},h_{2}(s)}
\end{displaymath}
It follows that there exists $f \in A$ such that for every
integer $m$
\begin{displaymath}
(\neg \psi)_{h,m} \equiv \bigwedge_{s=1}^{m} \neg
(\psi)_{h_{1}(s),h_{2}(s)}
\equiv \bigwedge_{s=1}^{m} \neg (\psi)_{f_{1}(s),f_{2}(s)}
\equiv (\neg \psi)_{f,m}
\end{displaymath}
Now, by definition of $O(n) \subseteq A$, there exists $g \in
O(n)$ such that $g\uparrow\{1,\ldots n\}=f\uparrow\{1,\ldots n\}$. We
obtain then
\begin{displaymath}
(\neg \psi)_{h,n} \equiv
\bigwedge_{s=1}^{n} \neg \mbox{ }(\psi)_{f_{1}(s),f_{2}(s)}
\equiv \bigwedge_{s=1}^{n} \neg\mbox{ }(\psi)_{g_{1}(s),g_{2}(s)}
\equiv
(\neg \psi)_{g,n}
\end{displaymath}
and this is the desired result.

{\bf Existential}. Fix a formula $\phi$ in
$L_{A}$ with free variables among the collection
\{$\vec v_{i}|i<\omega\}\cup \{\vec x\}$. Let $\vec
K=(K_{1},K_{2}\ldots )$ a corresponding vector of predicate symbols
with true sort. Consider the formula $
\psi(\vec x)=\exists \vec v \in \vec K \phi
( \vec v,\vec x)$. 

Recall:
$I(\psi)=I(\phi)$. Also, recall that $\forall n\in \omega\forall
h\in I(\psi(\vec x))$,
\begin{displaymath}
(\psi(\vec x))_{h,n}=\exists \vec v\in \vec K \phi_{h,n}(\vec
v,\vec x)
\end{displaymath}
Define then $Dense(\psi)=Dense(\phi)$. The verification of the desired property is trivial. \qed

A welcome consequence of this lemma is the following
corollary that states that the sets $I(\phi)$ are non empty.

\begin{cor}
For every formula $\phi \in
L_{A}$, $I(\phi)$ is non empty. \findef
\end{cor}
PROOF: By induction on formulas. The only interesting step is the
negation step. Recall that $I(\neg \phi)\subseteq (I(\phi)\times
\omega)^{\omega}$ is the collection of all maps $f=(f_{1},f_{2})$
such that:
\begin{equation}
\forall h \in I(\phi) \exists s \in \omega \mbox{, }
\phi_{h,f_{2}(s)} \equiv
\phi_{f_{1}(s),f_{2}(s)} \label{hh}
\end{equation}
 Lemma~\ref{setDense(phi)} states that $Dense(\phi)$ is countable. Hence
 the set $B=Dense(\phi) \times \omega \subseteq I(\phi) \times \omega$
 is countable. Let $f:\omega\longmapsto I(\phi)\times \omega$ be such
that $Im(f)=B$. It is easy to check using Lemma~\ref{setDense(phi)}
that $f$ verifies
statement~\ref{hh} above, hence $ f \in I(\neg \phi)$.
 \qed

\section{A Model Existence Theorem for $L_{A}$}
\label{comply}

Our intention is to prove  a model existence theorem for
uniform collections  of models using an ultraproduct construction very similar to the ultraproduct construction for  normed spaces (see for example \cite{Aksoy&Khamsi}).

In this section we fix a metric language $L(M_{\Re})$.

Let us recall some notation. Let $D$ a compact subset of 
the reals with the usual metric. 
Let $s=(s_{1},s_{2},\ldots s_{i}\ldots )$ $(i < \omega)$ be a sequence of
elements in $D$. It is well known that for every ultrafilter $U$
over $\omega$ there exists a unique $x\in D$ so that:
\begin{displaymath}
\forall \epsilon > 0 \mbox{ } \exists p \in U \mbox{ } \forall
i\in p
\mbox{  } |s_{i}-x| \leq \epsilon
\end{displaymath}
We will denote the element $x$ by $\lim_{i\in U} s_{i}$.

\begin{df} Ultraproduct Construction
\end{df}

Fix a uniform collection  of
models $\cal W$. Consider a sequence of models
$\{ \structi=(E_{i}, F^{i}_{\Re})
 | i < \omega \}$
  of $\cal W$. For every $i$, let $X_{i}$  be the universe of  the
 first order model $E_{i}$.
   For every $f$ in $\cal F$ we will denote by $f^{i}$
   the interpretation of this function symbol on the model
   $\structi$. Likewise for any relation symbol $C$ in $L$,  $C^{ i}$
 is the interpretation of $C$ in $\structi$.

   Recall that for every relation symbol $K \in L$, for
   every integer $n$, $K_{n}$ is the $(1/n) $ metric deformation of
   $K$.

Define  $X$ to be the set:
\begin{displaymath}
\{ g \in \prod_{i=1}^{\omega}X_{i} |  \mbox{ }
\exists K \in \mbox{$L$ a predicate such that }
\mbox{,  }
\forall n \in \omega\forall^{*} i \in \omega\mbox{, }
\structi \verify
K_{n}(g(i)) \}
\end{displaymath}
In other words, $X$ is the collection of all the sequences $g$ that
 eventually converge to some  predicate $K$.

Given a regular ultrafilter $\cal U $ over $\omega$, we define on
$X$ the usual equivalence relation:
For arbitrary $ x , y $ in $X$, $x \uequiv  y$ iff for
every $\epsilon>0 $ there exists $p \in \cal U$ such
that $\forall i \in p$ $\rho^{i}(x(i),y(i)) \leq \epsilon$. It is
easy to verify that this is truly an equivalence relation. We
will denote by $ [x]$ the equivalence class of $x$.

 The set $(X/  \uequiv)$ is endowed with a metric $\rho$ in the
 natural way:
 For every pair $ {[x]}, {[y]}$, find two representatives $ x,
  y$ such that there exist $K^{1},K^{2}$  predicates satisfying for
 every integer $i$: $\structi \verify  K^{1}(x(i))\wedge
   K^{2}
    (y(i))$
   (this is possible because $\cal U$ is regular). 

    By the property of directed sets of  uniform collections of models 
    we know
    that there exists a unary predicate $K$ such that for every
    integer $i$, $
\structi \verify K(x(i))\wedge K
 (y(i))$.
 Invoking again the properties of the uniform collections of models, we get
 that there exists a compact set $D$ in $\reals$ such that for
 every $i$, $
 \structi \verify \forall (z, v) \in (K,K)
 D(\rho^{i}(z, v))$. 
  It makes sense then to define
  \begin{displaymath}
  \rho( [x], [y])=\lim_{i\in \mbox{$\cal U$}} \rho^{i}
  ( x(i), y(i))
  \end{displaymath}
  It follows from the definition of a uniform collection of models that
  the metric function is well defined.

Using the set $(X/ \uequiv)$ we define the interpretations of
$L$:

\begin{itemize}
\item For every model valued function symbol $f$ in $L$ with arity
$a$ we define $f:(X/ \uequiv)^{a} \longmapsto (X/
\uequiv)$ by:
 \begin{quote}
  $\forall  \vec {[x]} \in (X/ \uequiv)^{a}$, $f(\vec {[x]})=
   [(f^{ 1}(\vec x(1)),\ldots f^{i}(\vec x(i)),\ldots )] $
   \end{quote}
   The properties of the uniform collection  of models guarantee that
   the image of $f$ is a subset of (X/ \uequiv) and that $f$ is well
   defined.

   \item For every real valued function symbol $f$ in $M_{\Re}$ with
 arity $a$
  we define $f:(X/ \uequiv)^{a}
   \longmapsto \reals$ by:
   \begin{quote}
    $\forall \vec {[x]} \in (X/ \uequiv)^{a}$,
     $f(\vec  {[x]})=lim_{i\in U} f^{ i}(\vec x(i) )$.
\end{quote}
The properties of the uniform collections  of models once again
guarantee that the images of $f$ are well defined.

\item For any relation symbol $C$ in $L$ with arity $a$ we define
 the interpretation of $C$ in $(X/
\uequiv)^{a}$ as follows:
\begin{quote}
$\vec {[x]} \in C$ iff for every integer $n$, $
\exists p \in \cal U $ such that for every  $i \in p$, $\structi
 \verify (\mbox{ }C(\vec x(i))\mbox{ })_n$
\end{quote}
\end{itemize}
Such a structure is denoted by $\ultrastructa$. \findef

The following remark follows directly from the previous
construction.

\begin{remark} The ultraproduct construction yields a model of $L(M_{\Re})$.
\label{model}\end{remark}
 
 Let $\cal {W}$ be a uniform collection  of
 models of $L$. Let $$ \{ \structi | i < \omega \}$$ 
     be a sequence of models of $\cal {W}$. Let $\cal U$ be a regular
     ultrafilter over $\omega$. Then the structure \ultrastructa  is
      a model of $L$.
       \findef

       PROOF:
       Let us verify first that \ultrastructa  is a model of $L$.

       Clearly $((X/\uequiv),d)$ is a metric space. Using the property
       of uniform continuity of the  uniform collection of models $\cal
       {W}$ we obtain that the functions $f$   are continuous.
It is easy to verify that the interpretation of the relation symbols
$C$ in $L$ with arity $a$ are closed subsets of
$(X/\uequiv)^{a}$ (in the product topology). 

Likewise, it is easy to verify that that the interpretations of the real valued functions 
symbols in $M_{\Re}$ are continuous functions.\qed

We introduce the following notation. Given any model
$\structa$, any $\vec a=(\vec a_{1},\ldots \vec a_{i},\ldots  )$ in
$\structa$ and any corresponding vector of true sort predicate
symbols $\vec K=(\vec K_{1},\ldots \vec K_{i},\ldots )$, by $\vec K(\vec
a)$ we understand $\bigwedge_{i=1}^{\infty} K_{i}(\vec a_{i})$.

The next lemma states the basic relationship between approximate 
formulas of $L_{A}$ and the ultraproduct.

\begin{lemma} Property of the Ultraproduct
\label{propultrapower}

 Fix  $\cal {W}$ a uniform collection  of models in 
$L(M_{\Re})$. Let $$ \{ \structi | i < \omega \}$$ 
    be  a sequence of models in $\cal {W}$. Fix $\cal U$ a regular 
    ultrafilter over $\omega$. For any formula $\phi(\vec x) \in
    L_{A}$, for any sequence $\{\vec a_{r}\}_{r=1}^{\infty}$ of
    vectors whose elements are in $(X/\uequiv) $, the following are
     equivalent:
      \begin{itemize}
      \item $\exists h \in I(\phi) \forall n \in \omega 
      \forall^{*} r\in \omega$
      $\ultrastructa \verify
      \phi_{h,n}(\vec a_{r})$
      \item
       $\exists h\in I(\phi)$ such that $\forall n \in \omega$ 
       $\forall^{*} r\in \omega$
$\exists p \in\cal U$ satisfying:
\begin{displaymath}
\forall i \in p \mbox{, } \structi  \verify \phi_{h,n}(\vec
a_{r}(i)). \mbox{ \findef}
\end{displaymath}
\end{itemize}
\end{lemma}

PROOF: By induction on formulas.

{\bf Atomic Formulas\/}. The proof follows directly from the definition of
 the
interpretation of the predicates in the model \ultrastructa, using
the facts that the predicates with fixed sort are closed, the
interpretation of the function symbols are continuous functions and that
 the number of branches in the tree of approximation of atomic formulas is 1.

{\bf Conjunction\/}. Both directions are direct.

{\bf Negation}. $\Rightarrow$. Suppose that $\exists g \in I(\neg
\phi)$  $\forall n\in \omega$  $\forall^{*} r \in \omega $ 
$\ultrastructa \verify (\neg
\phi(\vec a_{r}))_{g,n}$. By the definition of the 
negation step for approximate formulas  we get that for every $h$
in $I(\phi(\vec x))$ there exists an integer $n$ such that
$\forall^{*} r\in \omega$ $\ultrastructa \verify \neg \phi_{h,n}(\vec
a_{r})$.

 Using the induction hypothesis on the formula $\phi$, and the
 fact that $\cal U$ is an ultrafilter, we get that for all
 $h$ in $I(\phi)$ $\exists n$ $\forall^{*} r\in \omega$ there exists a $p$
 in $\cal U$ satisfying:
 \begin{displaymath}
 \forall i \in p \mbox{ }\structi \verify \neg \phi_{h,n}(\vec
 a_{r}(i))
 \end{displaymath}

 Recall that $Dense(\phi) \subset I(\phi)$ is a countable set of branches ``dense'' in $I(\phi)$ (see Lemma~\ref{setDense(phi)}). Define then the set $W=$
 \begin{displaymath}
 \{(h,n): h \in Dense(\phi) \wedge n \in \omega \wedge \mbox{
 }\forall^{*} r\in \omega \exists p \in \mbox{$\cal U$}\forall i\in p
 \mbox{, } \structi \verify \neg (\phi_{h,n}(\vec a_{r}(i))\}
 \end{displaymath}
Since the above set is countable, 
it is easy to obtain a
$g=(g_{1},g_{2}):\omega \longmapsto Dense(\phi) \times \omega$ such
that $Image(g)=W$. 
Using the properties of the dense countable set $Dense(\phi)$
(Lemma~\ref{setDense(phi)}) it is a standard procedure to verify that this
function is in $I(\neg\phi(\vec x))$.

Furthermore, by definition of $W$, $\forall n$ $\forall^{*} r\in \omega$
there exists a $p$ in $\cal U$ such that:
\begin{displaymath}
\structi \verify (\neg \phi(\vec a_{r}(i)))_{g,n} \equiv
\bigwedge_{s=1}^{n} \neg ( \phi(\vec a_{r}(i))
_{g_{1}(s),g_{2}(s)})
\end{displaymath}
This is the desired result.

$\Leftarrow$. Similar to the previous proof. Left to the reader.

{\bf Existential\/}. Here we use the full power of the induction
hypothesis on arbitrary sequences $\{a_{r}\}_{r=1}^{\infty}$.

$\Rightarrow$. Suppose that
$\exists h \in I(\phi) \forall n \forall^{*} r \in \omega$
\begin{displaymath}
\ultrastructa \verify  (\exists \vec x \in \vec K
\phi(\vec a_{r},\vec x))_{h,n}
\end{displaymath}

By the definition of approximate formulas, the above statement
implies that $\exists h \in I(\phi) \forall n \forall^{*} r\in \omega $
\begin{displaymath}
\ultrastructa \verify \exists \vec x \in \vec K
(\phi(\vec a_{r},\vec x))_{h,n}
\end{displaymath}
It follows that for every integer $r$ there exists a vector
\begin{displaymath}
\vec d_{r}
\end{displaymath}
such that $\forall n \forall ^{*} r$
\begin{displaymath}
\ultrastructa \verify \vec K(\vec d_{r}) \wedge \phi_{h,n}(\vec
a_{r},\vec d_{r})
\end{displaymath}
Furthermore, by the definition of $\sim$ we can select the elements of 
$\vec d_{r}$ in such a way that for every integer $i$, $\structi \verify 
\vec K(\vec d_{r}(i))$.

Using the induction hypothesis, we obtain then that there exists
$h \in I(\phi)$ such that $\forall n$ $\forall^{*} r$  there
exists $p \in \cal {U}$ with :
\begin{displaymath}
\forall i \in p \mbox{ } \structi \verify  (\vec K)
(\vec d_{r}(i))
\wedge \phi_{h,n}(\vec a_{r}(i),\vec d_{r}(i))
\end{displaymath}
But this implies that $\forall n$, $\forall^{*}r$, $\exists p \in \cal U$,
\begin{displaymath}
\forall i \in p \mbox{ }\structi \verify
 \exists \vec x\in
 \vec
 K \phi_{h,n}(\vec a_{r}(i),\vec x)
 \end{displaymath}
 This is the desired result.

 $\Leftarrow$: Similar to the previous proof. Left to the reader.

This completes the proof of the lemma. \qed

We use the previous result to show that the ultraproduct is a
rich model, i.e. a model where $\apptruth$ and $\verify$ coincide
for all formulas in $L_{A}$.

\begin{theo} The Ultraproduct is Rich for $L_{A}$
\label{rich}

 Fix  a uniform collection  of models
 $\cal W$ in $L(M_{\Re})$. Let $\{ \structi\}_{i=1}^{\omega}$
 be a sequence of
 models in $\cal {W}$. Fix $\cal U$ a regular
 ultrafilter over $\omega$. Then for every formula $\phi(\vec x)$
 in $L_{A}$ and for every $\vec a$, a vector of elements in $\ultrastructa$   with finite or infinite arity, the following are equivalent:
 \begin{itemize}
 \item
 \ultrastructa \apptruth $ \phi(\vec a)$
\item \ultrastructa \verify $\phi(\vec a)$ \mbox{ \findef}
\end{itemize}
\end{theo}

PROOF: By induction in formulas. The atomic, negation and countable 
conjunction cases are direct and we leave them to the reader. 

For the existential case, there is only one interesting
direction. Suppose that $\ultrastructa \apptruth \exists \vec v \in \vec K \phi(\vec v,\vec a)$.

 From Lemma~\ref{propultrapower} it follows that there exists
  $h \in I(\phi)$ such that for every integer $n$  $ \exists p \in
  \cal U$ such that:
  \begin{displaymath}
  \forall i \in p \mbox{ }
  \structi \verify \exists\vec v \in \vec K \phi_{h,n}(\vec v,\vec
  a(i))
  \end{displaymath}

  Since the ultrafilter is regular, it is easy to
  construct by diagonalization a sequence $\vec b$ of
  vectors of functions (elements of \ultrastructa) corresponding to
  $\vec v$ and such that $\forall n \in \omega$ $\exists
   p \in \cal {U}$ $\forall i \in p$ :
   \begin{displaymath}
\structi \verify \vec K(\vec b(i))\wedge
\phi_{h,n}(\vec b(i),\vec a(i)))
\end{displaymath}

 We again invoke Lemma~\ref{propultrapower}
  to obtain that $
  \ultrastructa \apptruth \vec K(\vec b) \wedge \phi(\vec b,\vec a) $.

  By induction hypothesis we get $
  \ultrastructa \verify
   \vec K(\vec b) \wedge \phi(\vec b,\vec a)$,
and this implies $
   \ultrastructa \verify \exists \vec v\in \vec K \phi(\vec
   v,\vec a)$. 
This completes the proof.
\qed

The main consequence of the previous lemmas is the following
Model Existence Theorem.

 \begin{theo} Model Existence Theorem for $L_{A}$
\label{compact}

 Let $\cal
M$ be a uniform collection of models for $L(M_{\Re})$. Fix a sentence 
$\phi \in L_{A}$. Suppose that there exists a branch 
$h \in I(\phi)$ such that for every integer $n$ there exists 
$\structn \in \cal {M}$ such that:
\begin{displaymath}
\structn \verify \phi_{h,n}
\end{displaymath}

Then there exists a  model $\structa $ in $\cal {M}$
 such that:
\begin{displaymath}
\structa \apptruth \phi  \mbox{ and } \structa \verify \phi
\mbox{ \findef}
\end{displaymath}
\end{theo}

PROOF: Select a regular ultrafilter over $\omega$
 and
 apply Lemma~\ref{propultrapower} to the
  sequence of models $\{\structn\}_{n=1}^{\infty}$. We obtain that $
  \ultrastructa \apptruth \phi $.
 
Applying now Theorem~\ref{rich} to $\ultrastructa$ we get that 
   $ \ultrastructa  \verify \phi $.

It remains to show that \ultrastructa is in the uniform
collection $\cal {W}$. It is enough 
to check that \ultrastructa satisfies the same
uniform bounds for the function symbols in $L(M_{\Re})$, the same
inclusion relationship for the true sort predicates and the same
uniform continuity property for the function symbols that the
models in $\cal {W}$. Left to the reader.\qed

\begin{remark} $\empty$
\end{remark}

Recall that in Remark~\ref{logichenson} we mentioned that for formulas
 $\phi \in L_{PB}$ the tree of approximations has only one branch $h$  (up to logical equivalence), and 
that for this branch, the approximate formulas $\{\phi_{h,n}\}_{n=1}^{\infty} $ coincide with Henson's approximate 
formulas $\{\phi_{n}\}$. 

It follows then that we have, as a corollary of the Model Existence
 Theorem, the following version of Henson's Compactness Theorem:

\begin{cor} Henson's Compactness Theorem
\end{cor}

Fix $L(\{\rho, ||.||\}\cup M_{\Re})$ with $L$ and $\{\rho,||.||\}$ as
 in Example~\ref{example1}. Let $\cal M$ a uniform collection of normed
 structures in the above language. Fix a sentence $\phi \in L_{PB}$.
 If for every integer $n$ there exists a model $\structn \in \cal M$ 
with $\structn \verify \phi_{n}$, then there exists a model $\structa
 \in \cal M$, such that $\structa \apptruth \phi$ and $\structa \verify
 \phi$.\findef

\bigskip

Beaver College, Glenside, PA,

ortiz@beaver.edu
\medskip

Mathematical Subject Classification, 1991: 

 Primary 03C65

 Secondary  46B08, 46B20
\medskip

Keywords: approximate truth, compactness theorem, normed space structures.

\end{document}